\documentclass[12pt]{amsart}
\usepackage{amssymb,amsmath, amsthm,latexsym}
\usepackage{a4wide}
\newcommand{\cal}[1]{\mathcal{#1}}
\theoremstyle{plain}

\newtheorem{lem}{Lemma}[section]
\newtheorem{theo}[lem]{Theorem}
\newtheorem{prop}[lem]{Proposition}
\newtheorem{corollary}[lem]{Corollary}

\parskip=\medskipamount
\font\k=cmr7

  \newcommand {\cu}{\mbox{\k cusp}}
  \newcommand {\di}{\mbox{\k disc}}

  \newcommand {\fin}{\mbox{\k fin}}
  \newcommand {\Gr}{G(\R)}
  \newcommand {\C}{{\mathbb C}}
  \newcommand {\N}{{\mathbb N}}
  \newcommand {\R}{{\mathbb R}}
  \newcommand {\Z}{{\mathbb Z}}
  \newcommand {\Q}{{\mathbb Q}}
  \newcommand {\A}{{\mathbb A}}
  
  \newcommand {\gf}{{\mathfrak g}}

  \newcommand {\ho}{{\mathfrak o}}

 \newcommand {\mX}{{\mathfrak X}}

 \newcommand {\mG}{{\mathfrak G}}
  
\renewcommand {\H}{{\mathcal H}}
  
  \newcommand {\Co}{{\mathcal C}}
  \newcommand {\cO}{{\mathcal O}}

\newcommand {\ba}{\backslash}

\newcommand{\Tr}{\operatorname{Tr}}
\newcommand{\End}{\operatorname{End}}

\newcommand{\tr}{\operatorname{tr}}
\newcommand{\Id}{\operatorname{Id}}
\newcommand{\Hom}{\operatorname{Hom}}

\newcommand{\vol}{\operatorname{vol}}

\newcommand{\SL}{\operatorname{SL}}

\newcommand{\SO}{\operatorname{SO}}

\newcommand{\PSL}{\operatorname{PSL}}

\begin{document}
\title[Weak Weyl's Law for congruence subgroups]{\Large\bf Weak Weyl's Law for congruence subgroups}
\date{\today}

\author{Jean-Pierre Labesse and Werner M\"uller}
\address{Institut de Math\'ematiques de Luminy, Universit\'e Aix-Marseille II,
Av. de Luminy, Case Postale 907, 13288 Marseille, CEDEX 9, France}
\email{labesse@math.jussieu.fr}
\address{Universit\"at Bonn\\
Mathematisches Institut\\
Beringstrasse 1\\
D -- 53115 Bonn, Germany}
\email{mueller@math.uni-bonn.de}

\maketitle

\section[Introduction]{Introduction}
\setcounter{equation}{0}

Let $G$ be a connected and simply connected semisimple algebraic group defined
over $\Q$ and let $Z_G$ be its center. Let $\Gamma$ be an arithmetic subgroup 
of $\Gr$. For simplicity we assume in the introduction that $\Gamma$ is 
torsion free. We denote by $K_\infty$ a maximal compact subgroup of $\Gr$.
We endow $\Gr$ and $K_\infty$ with Haar measures so that the symmetric space
$X=\Gr/K_\infty$ has the canonical measure, i.e., the measure defined using
 the Killing form.  Let 
$$H^\Gamma_{\cu}=L^2_{\cu}(\Gamma\ba \Gr)$$
denote  the closure in $L^2(\Gamma\ba \Gr)$ of the subspace 
spanned by all $\Gamma$-cuspidal automorphic forms. 
We denote by $\rho_\infty$ the right regular
representation of $\Gr$ in $H^\Gamma_{\cu}$. Let $(\sigma,V_\sigma)$ be an
 irreducible unitary representation of $K_\infty$. Set
$$H_{\cu}^\Gamma(\sigma)=(H_{\cu}^\Gamma\otimes V_\sigma)^{K_\infty}.$$
Let $\Omega$ be the Casimir operator for $\Gr$. Then 
$\rho_\infty(\Omega)\otimes \Id$ induces a selfadjoint operator in the Hilbert
space $H^\Gamma_{\cu}(\sigma)$ which has pure point spectrum.  We denote by 
$\Lambda_{\cu}(\sigma)$ the set of eigenvalues of 
$\rho_\infty(\Omega)\otimes \Id$ in $H^\Gamma_{\cu}(\sigma)$. 
Let $N^\Gamma_{\cu}(T,\sigma)$ be the counting 
function of the cuspidal spectrum $\Lambda_{\cu}(\sigma)$, i.e.,  
$N^\Gamma_{\cu}(T,\sigma)$ is the 
number of eigenvalues $\lambda$ of $\rho_\infty(\Omega)\otimes\Id$ in 
$H^\Gamma_{\cu}(\sigma)$, counted with their multiplicities, 
with $|\lambda|\le T$. 
We shall also consider the discrete subspace 
$$H^\Gamma_{\di}=L^2_{\di}(\Gamma\ba \Gr),$$
and we define $H^\Gamma_{\di}(\sigma)$ and the corresponding eigenvalue 
counting function $N^\Gamma_{\di}(T,\sigma)$ similarly.
It is an important question in the theory of automorphic forms to understand
 the asymptotic behavior of the counting function as $T\to\infty$. Let
$d$ be the dimension of the symmetric space $X=\Gr/K_\infty$ and let 
${\bf \Gamma}$ denote the Gamma function. Let
$$c_\sigma(\Gamma)=
\frac{\dim(\sigma)\vol(\Gamma\ba X)}{(4\pi)^{d/2}{\bf \Gamma}(d/2+1)}$$
be Weyl's constant. A conjecture of Sarnak \cite{Sa} 
states that the counting function of
the cuspidal spectrum satisfies Weyl's law, i.e., 
\begin{equation}\label{1.1}
\lim_{T\to\infty}\frac{N^\Gamma_{\cu}(T,\sigma)}{T^{d/2}}=\lim_{T\to\infty}
\frac{N^\Gamma_{\di}(T,\sigma)}{T^{d/2}}=c_\sigma(\Gamma).
\end{equation}

This conjecture has been established in some special cases. First of all,
it was Selberg \cite{Se} who proved it for congruence subgroups of
$\SL(2,\Z)$ and $\sigma=1$. Other cases for which the conjecture has been
established are Hilber modular groups \cite{Ef},
congruence subgroups of $\SO(n,1)$ \cite{Rez}, $\SL(3,\Z)$ \cite{Mil}, and
in particular, the conjecture  was proved  in \cite{Mu2} for  principal 
congruence subgroups of $\SL(n,\Z)$ and arbitrary $\sigma$. 

A possible approach to prove (\ref{1.1}) is through the Selberg trace formula
applied to a function constructed from the heat kernel of the corresponding
Bochner-Laplace operator. Let $E_\sigma\to 
X$ be the  homogeneous vector bundle associated to $\sigma$
and let $\Delta_\sigma$ be the  elliptic differential operator induced by
$-R(\Omega)\otimes\Id$ in $C^\infty(X,E_\sigma)$, where $R$ denotes the right
regular representation of $\Gr$ in $C^\infty(\Gr)$. Using the heat 
kernel of $\Delta_\sigma$,  one 
produces a 1-parameter family of functions $h_t$ on $\Gr$, depending on $t>0$,
 such that for any irreducible unitary representation $(\pi,\H(\pi))$ of 
$\Gr$ one has
\begin{equation}\label{1.2}
\Tr\pi(h_t)=e^{t\pi(\Omega)}\dim\Hom_{K_\infty}(\H(\pi),V_\sigma^*).
\end{equation}
Given $\lambda\in\Lambda_{\cu}(\sigma)$, let $m(\lambda)$ denote its
 multiplicity. Then it follows from (\ref{1.2}) that
$$\Tr\left(\rho_\infty(h_t)|H^\Gamma_{\cu}\right)
=\sum_{\lambda\in\Lambda_{\cu}(\sigma)}m(\lambda)e^{t\lambda}.$$
The idea is now to study the asymptotic behavior of the left hand side as
$t\to0$. 
Using the properties of the heat kernel, it follows that 
\begin{equation}\label{1.3}
h_t(1)\sim  \frac{\dim(\sigma)}{(4\pi t)^{d/2}\vol(K_\infty)}
\end{equation}
as  $t\to 0$. A good control of the
 trace formula should allow to prove that as $t\to0$, the constant term in
the asymptotic expansion  of 
$$t\longmapsto t^{d/2}\Tr\left(\rho(h_t)|H^\Gamma_{\cu}\right)$$
 is the contribution of the unit element to the trace formula:
$$t^{d/2}\Tr\left(\rho(h_t)|H^\Gamma_{\cu}\right)=t^{d/2}\vol(\Gamma\ba
 G(\R))h_t(1)+o(1)$$
(we recall that we have assumed that $\Gamma$ is torsion free). 
The desired estimate is then an immediate consequence of (\ref{1.3}) combined
with the Tauberian theorem. The second named author has been able to apply successfully this method for 
$\Gamma$ a principal congruence subgroup for $G=\SL(n)$. 
However,  many technical difficulties are still to overcome in order to prove
(\ref{1.1}) for a general $G$.  

The purpose of this paper is to prove a weaker result which holds for every
$G$. Recall that an upper bound, which holds for arbitrary $G$ and $\Gamma$, 
is already known thanks to Donnelly \cite{Do}: 
$$\limsup_{T\to\infty}\frac{N^\Gamma_{\cu}(T,\sigma)}{T^{d/2}}\le 
c_\sigma(\Gamma).$$
By working with a simple form of the trace formula, we shall get, for a general
$G$, a lower bound that depends on the choice of a set $S$ of primes containing
at least two finite primes. For every such set $S$  we shall define a certain 
constant $c_S(\Gamma)\le 1$, which is non zero for $\Gamma$ deep enough. Let 
$Z_\Gamma=\Gamma\cap Z_G(\R)$. Note that $Z_\Gamma\subset K_\infty$.
Then our main result is the following theorem.

\begin{theo}\label{th1.1} Let $G$ be an almost simple  connected and simply 
connected semisimple algebraic group defined over $\Q$ such that $\Gr$
 is non compact. Then for every $\Gamma$ 
and every $\sigma$ such that $\sigma|_{Z_\Gamma}=\Id$ we have
$$c_\sigma(\Gamma)c_S(\Gamma)\le \liminf_{T\to\infty}
\frac{N^\Gamma_{\cu}(T,\sigma)}{T^{d/2}}.$$
\end{theo}

This will turn out to be very easy to prove but the price to pay for 
simplicity is that this lower bound is far from being sharp; in particular 
the lower
bound is non trivial, i.e. $c_S(\Gamma)>0,$ only for discrete subgroups 
$\Gamma$ that are deep enough (see below for a more precise statement). 
We note that a
weaker form of this result is due to Piatetski-Shapiro \cite{PS} and for 
$G=\SL(2,\R)$ also to Venkov \cite{Ve}.
In \cite{PS} it is proved that for
every $\Gamma$ there exists a subgroup of finite index $\Gamma'$ such
that the space of Maass cusp forms with respect to $\Gamma'$ is infinite
dimensional. 

Theorem \ref{th1.1} is a consequence of Theorem \ref{th4.2}. 
To get the lower bound we shall prove a sharp estimate for a certain strict
subspace $H^{\Gamma,S}_{\cu}$ of $H^\Gamma_{\cu}$ that will
be defined via the ad\`elic picture. From now on, for simplicity of exposition,
let us assume that $G$ is absolutely almost simple and that $G(\R)$ is non 
compact. We consider a discrete subgroup $\Gamma$ such that
$$\Gamma=K_{\fin}\cap G(\Q),$$
where
$$K_{\fin}=\prod K_p$$
is a decomposable open compact subgroup of the group of finite ad\`eles
$G(\A_{\fin})$. Let $S$ be a finite set of primes of 
cardinality at least $2$ and let
$$L^2_{\cu}(G(\Q)\ba G((\A),S)\subset L^2_{\cu}(G(\Q)\ba G(\A))$$
be the $G(\A)$-module generated by cusp forms that are orthogonal to the 
trivial representation of $G(\A)$ and on which $G_S=G(\Q_S)$ acts by the 
Steinberg representation. Then we define 
$$H^{\Gamma,S}_{\cu}=L^2_{\cu}(G(\Q)\ba G(\A),S)^{K_{\fin}}$$
and
$$H_{\cu}^\Gamma(\sigma,S)=\left(L^2_{\cu}(G(\Q)\ba G(\A),S)\otimes 
V_\sigma\right)^{K}=\left(H_{\cu}^{\Gamma,S}\otimes
V_\sigma\right)^{K_\infty}.$$
Let $\Lambda_{\cu}(\sigma,S)$ be the spectrum of $\rho_\infty(\Omega)\otimes
\Id$ in $H^\Gamma_{\cu}(\sigma,S)$. It consists of a discrete set of 
eigenvalues. Given $\lambda\in\Lambda_{\cu}(\sigma,S)$  denote by 
$m(\lambda)$ the multiplicity of the eigenvalue $\lambda$. Then it follows
from  (\ref{1.2}) that
$$\Tr(\rho_\infty(h_t)|H^{\Gamma,S}_{\cu})
=\sum_{\lambda\in\Lambda_{\cu}(\sigma,S)}m(\lambda)e^{t\lambda}.$$ 
Using the trace formula, we shall prove in 4.2 that
$$\Tr(\rho_\infty(h_t)|H^{\Gamma,S}_{\cu})\sim c_S(\Gamma)
\frac{\dim\sigma\vol(\Gamma\ba X)}{(4\pi t)^{d/2}}$$
as $t\to0$. For $T\ge0$ 
let $\Lambda_{\cu}(\sigma,S;T)$ be the subset of all $\lambda\in
\Lambda_{\cu}(\sigma,S)$ with $|\lambda|\le T$. Set
$$N^\Gamma_{\cu}(T,\sigma,S)
=\sum_{\lambda\in\Lambda_{\cu}(\sigma,S;T)}m(\lambda).$$
If $c_S(\Gamma)\neq 0$ the Tauberian theorem yields the estimate:
$$\lim_{T\to\infty}\frac{N^\Gamma_{\cu}(T,\sigma,S)}{T^{d/2}}=
c_\sigma(\Gamma)c_S(\Gamma).$$
This is our Theorem 4.2. Now clearly $N^\Gamma_{\cu}(T,\sigma,S)\le 
N^\Gamma_{\cu}(T,\sigma)$ and hence 
$$c_\sigma(\Gamma)c_S(\Gamma)\le\liminf_{T\to\infty} 
\frac{N^\Gamma_{\cu}(T,\sigma)}{T^{d/2}}.$$
We still have to observe that $c_S(\Gamma)\neq0$ for subgroups $\Gamma$
such that  at primes $p\in S$ the open compact subgroup $K_p$ is 
small enough, i.e., if it is a subgroup of some minimal parahoric subgroup 
$I_p$ of $G_p=G(\Q_p)$.

\section[A simple trace formula]{A simple trace formula}
\setcounter{equation}{0}

From now on we shall work in the ad\`elic picture. We choose a Haar measure,
 for example the Tamagawa measure, on $G(\A)$ which is given as the product of
measures on each $G_p=G(\Q_p)$. We fix a decomposable open compact subgroup
$K_{\fin}$ of $G(\A_{\fin})$. Set
$$K=K_\infty K_{\fin}.$$
Then $K$ is a compact subgroup of $G(\A)$. 
Let $f\in C^\infty_c(G(\A))$ be $K$-finite. The Arthur-Selberg trace formula
 in its invariant form is the equality
$$\sum_{\ho\in\cO}I_{\ho}(f)=I(f)=\sum_{\chi\in\mX}I_\chi(f)$$
of a geometric and a spectral expansion for an invariant linear
form $f\longmapsto I(f)$ (see \cite{Ar}). 

Let  $\Sigma$ be some large enough set of primes, 
containing in particular the real place and all ramified
primes. 
We shall consider functions $f$ that are decomposable 
$$f=\otimes f_p$$ 
and such that for ${p\notin\Sigma}$, $f_p$ is the characteristic function of
 a hyperspecial maximal compact subgroup  divided by its volume.
We put, as usual, 
$$f^\Sigma=\mathop{\otimes}_{p\notin\Sigma}f_p$$
and  
$$f_\Sigma=\mathop{\otimes}_{p\in\Sigma}f_p$$
so that $f=f_\Sigma\otimes f^\Sigma$.
The main observation is the following result of Arthur:
when applied to a decomposable function such that $f_p$ is 
cuspidal, in the sense of Arthur, at two primes $p\in S$ 
the trace formula is ``simple". 

\begin{prop}\label{p2.1}
If $f$ is cuspidal at two places, the linear form
$I(f)$  has a discrete spectral expansion:
$$I(f)=\sum a(\pi)\tr\pi(f)$$
where the sum is over certain subrepresentations of 
parabolically induced representations
from representations in the discrete spectrum of Levi subgroups. In
particular, if $\pi$ occurs in the discrete spectrum for $G(\A)$ then
$a(\pi)=m(\pi),$
where $m(\pi)$ denotes the multiplicity with which $\pi$ occurs in the 
discrete spectrum.
The geometric expansion is given by a linear combination of
ordinary orbital integrals (over $G_\Sigma$)
$$I(f)=\sum_{\gamma\in\mG_\Sigma} a(\gamma,\Sigma)
{\cO}_\gamma(f_\Sigma),$$
where $\mG_\Sigma$ is a set of representatives of conjugacy classes
with elliptic semisimple part, ${\cO}_\gamma(f_\Sigma)$
is the orbital integral of $f_\Sigma$ and $a(\gamma,\Sigma)$
is some scalar. If $\gamma$ is  elliptic we have
$$ a(\gamma,\Sigma){\cal O}_\gamma(f_\Sigma)= a(\gamma){\cO}_\gamma(f),$$
where the scalar $a(\gamma)$ is independent of $\Sigma$~.
Moreover if the support of $f_\Sigma$ remains in a fixed compact set
then, in the  sum over $\mG_\Sigma$~, only a finite number of elliptic
conjugacy classes in $G(\Q)$ may give a non trivial contribution.
\end{prop}

\begin{proof}
This is the content of \cite[Corollary 7.2]{Ar}. We omit
the sum over $t\ge1$ that occurs in \cite{Ar} since the sum is now
known to be absolutely convergent by \cite{Mu1}.
\end{proof}

We shall denote by $1(f)$ the operator defined by $f$ in the trivial
 representation; this is a scalar operator in a one dimensional space:
$$1(f)=\int_{G(\A)}f(x)\,dx.$$
Let $p$ be a finite place. The existence
of (normalized) pseudo-coefficients of the Steinberg representation is
a particular case of a general result due to Kazdan, 
but, up to a sign, they can also be viewed
as Euler-Poincar\'e functions at non archimedean places
and an explicit construction has been given by Kottwitz
\cite{Ko}. A pseudo-coefficient of the Steinberg representation
is a $K_p$-finite function $f_p\in C^\infty_c(G_p)$ such that 
$$\Tr \pi_p(f_p)=0$$ 
for all irreducible unitary representations $\pi_p$
unless either $\pi_p$ is the Steinberg representation of $G_p$
in which case 
$$\Tr \pi_p(f_p)=1$$
or $\pi_p=1_p$  the trivial representation in which case 
$$1_p(f_p)=(-1)^q$$
for some integer $q$.
Moreover, \cite[Theorem 2]{Ko} shows that orbital integrals 
${\cO}_\gamma(f_p)$ of non elliptic elements $\gamma$ vanish.

Let $S$ be a finite set of finite primes. Let 
$$L^2_{\cu}(G(\Q)\ba G(\A),S)\subset L^2_{\cu}(G(\Q)\ba G(\A))$$
be the $G(\A)$-module generated by cusp forms
that are orthogonal to the trivial
representation of $G(\A)$ and on which  
$$G_S:=\prod_{p\in S}G(\Q_p)$$ 
acts by the Steinberg representation. Let $\rho$ denote the regular 
representation of $G(\A)$ in $L^2_{\cu}(G(\Q)\ba G(\A))$.

\begin{prop}\label{p2.2}
 Let $G$ be a simply connected almost simple algebraic group 
defined over $\Q$. Assume that $S$ is a set of finite primes
containing at least two primes. Consider
$$f=f_\infty\otimes f_{\fin,S}\otimes f_S,$$
where $f_\infty\in C^\infty_c(G_\infty)$, 
$f_{\fin,S}$ is a locally constant compactly supported function on  
$G(\A_{\fin}^S)$, the group of finite ad\`eles outside $S$,  and $f_S$ is a 
pseudo-coefficient of the Steinberg representation
of $G_S$. Then  one has
$$\sum_{\gamma\in\mG_e} a(\gamma){\cal O}_\gamma(f)=
\Tr(\rho(f)|L^2_{cusp}(G(\Q)\ba G(\A),S))+1(f),$$
where $\mG_e$ is a set of representatives of elliptic conjugacy classes
in $G(\Q)$. 
\end{prop}

\begin{proof} 
By a standard argument, automorphic representations in the discrete spectrum 
with a Steinberg component at some place $p$ must be cuspidal. In fact
Langlands description of the discrete spectrum implies that representations
that are in the non cuspidal discrete spectrum 
arise from residues of Eisenstein series for non unitary parameters 
and hence are non tempered
at all places, while the Steinberg representation is tempered.
Moreover, automorphic representations with a
Steinberg component at some place $p$ cannot be
subrepresentations of parabolically induced representations 
from representations in the discrete spectrum of Levi subgroups
for a proper parabolic subgroup. Suppose that $G_p$ is non compact.
Since our group $G$ is absolutely almost simple and 
 simply connected, the strong approximation theorem holds and 
hence  $G(\Q)G_p$ is dense in $G(\A)$. Therefore an automorphic
representations that is trivial at $p$ is necessarily trivial. 
On the other hand, if the group $G_p$ is  compact
then all automorphic representations are cuspidal.
Hence an automorphic representations that is trivial at one place is
either trivial or cuspidal. 
We now use Proposition \ref{p2.1} on the simple form of the trace
 formula with $\Sigma$ large enough, containing $S$. The above
discussion shows that the spectral expansion of $I(f)$ 
is given by the sum of the trace in the cuspidal spectrum
plus the trace in the trivial representation:
$$I(f)=\Tr(\rho(f)|L^2_{cusp}(G(\Q)\ba G(\A),S))+1(f).$$
The geometric expansion of  
$I(f)$ is a linear combination of orbital integrals 
(over $G_\Sigma$) and, as recalled above,
\cite[Theorem 2]{Ko} shows that orbital integrals
of non elliptic elements vanish at $p\in S$. Hence
$$I(f)=\sum_{\gamma\in\mG_e} a(\gamma){\cO}_\gamma(f),$$
where $\mG_e$ is a set of representatives of elliptic conjugacy classes
in $G(\Q)$.
\end{proof}

Let 
$$H^{\Gamma,S}_{\cu}=L^2_{cusp}(G(\Q)\ba G(\A),S)^{K_{\fin}}$$ 
be the Hilbert subspace of
 $L^2_{\cu}(G(\Q)\ba G(\A))$ generated by vectors of automorphic 
representations that are Steinberg at $p\in S$
and that are $K_{\fin}$-invariant.
Let $e_{\fin}$ be the characteristic function
of $K_{\fin}$, our chosen open compact subgroup at finite ad\`eles
  divided by its volume:
$$e_{\fin}=\frac{1}{\vol(K_{\fin})}\chi_{K_{\fin}}.$$
Let
$$K_S=\prod_{p\in S}K_p,\quad K_{\fin,S}=\prod_{p\notin S\cup\{\infty\}}K_p.$$
Let $e_S$ (resp. $e_{\fin,S}$) be the characteristic function
of $K_S$ (resp. $K_{\fin,S}$) divided by its volume.
Then 
$$e_{\fin}=e_S\otimes e_{\fin,S}.$$

\begin{corollary}\label{c2.3}
Assume that $f$ is of the form
$$f=f_\infty\otimes f_S\otimes e_{\fin,S}$$
where 
$f_S$ is a pseudo-coefficient of the Steinberg representation
of $G_S$.
Let  $c(K_S)$ be the dimension of  the space of $K_S$ fixed vectors in the 
Steinberg representation of $G_S$ and let
$$b(f_\infty,S)=\Tr 1(f_\infty\otimes f_S),$$ 
where $1$ is the trivial
representation of $G_\infty\times G_S$. Then  one has
$$
\Tr\left(\rho_\infty(f_\infty)|H^{\Gamma,S}_{\cu}\right)=
c(K_S)\left(\sum_{\gamma\in\mG_e} a(\gamma){\cO}_\gamma(f)
-b(f_\infty,S)\right).$$
\end{corollary}
\begin{proof} We first observe that for any admissible representation
 $\tau$ of $G(\A_{\fin})$, the operator
$\tau(e_{\fin})$ is the projection on the subspace of $K_{\fin}$-invariant 
vectors. Hence,
$$\Tr(\rho_\infty(f_\infty)|H^{\Gamma,S}_{\cu})=
\Tr(\rho_\infty(f_\infty\otimes e_{\fin})|L^2_{\cu}(G(\Q)\ba G(\A),S)).$$
However, if $\pi_S$ is the Steinberg representation of $G_S$ then 
$$\Tr \pi_S(e_S)=c(K_S)\Tr\pi_S(f_S)$$
and hence we get
$$\Tr(\rho_\infty(f_\infty)|H^{\Gamma,S}_{\cu})=
c(K_S)\Tr(\rho(f)|L^2_{\cu}(G(\Q)\ba G(\A),S)).$$
The assertion is then an immediate consequence of Proposition \ref{p2.2}.
\end{proof}

\section[The heat kernel]{The heat kernel}
\setcounter{equation}{0}

In this section we establish a number of facts about the  heat kernel of 
a Bochner-Laplace operator on the Riemannian symmetric space 
$X=\Gr/K_\infty$. 
Let $E_\sigma\to X$ be the homogeneous vector bundle associated to $\sigma$. 
Let $\Omega\in{\mathcal Z}({\gf_\C})$ be the Casimir element of $\Gr$
and let $R$  be the right regular representation of $\Gr$ on
$C^\infty(\Gr)$.
Denote by $\Delta_\sigma$ the second order elliptic 
differential  operator induced by $-R(\Omega)\otimes{\Id}$ in 
$C^\infty(X,E_\sigma)\cong(C^\infty(\Gr\otimes V_\sigma)^{K_\infty}$. If
$\nabla^\sigma$ denotes the canonical invariant connection on $E_\sigma$, then
$$\Delta_\sigma=(\nabla^\sigma)^*\nabla^\sigma-
\lambda_\sigma{\Id},$$
where $\lambda_\sigma=\sigma(\Omega_{K_\infty})$  is the Casimir eigenvalue of
$\sigma$. Hence 
$\Delta_\sigma :C^\infty_c(X,E_\sigma)\to L^2(X,E_\sigma)$ is
essentially selfadjoint and bounded from below. We continue to denote its
 unique selfadjoint extension by $\Delta_\sigma$. Let
$\exp(-t\Delta_\sigma)$, $t\ge 0$, be the associated heat semigroup.
The heat operator is a smoothing operator on $L^2(X,E_\sigma)$ which commutes
 with the representation of $\Gr$ on 
$L^2(X,E_\sigma)$.
Therefore, it is of the form
$$
(e^{-t\Delta_\sigma}\varphi )(g)=
\int_{\Gr}H_t(g^{-1}g_1)(\varphi(g_1))dg_1, \quad
g\in \Gr,
$$
where $\varphi\in (L^2(\Gr)\otimes V_\sigma)^{K_\infty}$ and
$H_t\colon \Gr\to {\End}(V_\sigma)$ is in $L^2\cap C^\infty$ and 
satisfies the covariance property
$$
H_t(g)=\sigma(k)H_t(k^{-1}gk')\sigma(k')^{-1},\quad
{\mbox for}\;g\in \Gr,\;k,k'\in K_\infty.
$$
Let ${\mathcal C}^1(\Gr)$ be Harish-Chandra's space of
integrable rapidly decreasing functions on $\Gr$. Then 
\begin{equation}\label{3.1}
H_t\in({\mathcal C}^1(\Gr)\otimes{\End}(V_\sigma))^{K_\infty\times
K_\infty}
\end{equation}
\cite[Proposition 2.4]{BM}.

Next we need to estimate the covariant derivatives of the heat kernel. 
Let $\nabla$ denote the Levi-Civita connection on $\Gr$ with respect to
the invariant metric on $\Gr$. Let $d(x,y)$ be the geodesic distance of 
$x,y\in X$. Define the function $r(g)$ on $\Gr$ by
\begin{equation}\label{3.2}
r(g)=d(gK_\infty,K_\infty),\quad g\in G(\R).
\end{equation}
The following proposition is proved in \cite[Proposition 2.1]{Mu2}.
\begin{prop}\label{p3.1}
Let $a=\dim \Gr$,  $l\in \N_0$ and $T>0$. There exist $C,c>0$ such that
$$\parallel\nabla^l H_t(g)\parallel\le C t^{-(a+l)/2}
\exp\left( -\frac{cr^2(g)}{t}\right)$$
for all $0<t\le T$ and $g\in \Gr$. 
\end{prop}
Define the one-parameter family of smooth functions $h_t$
on $\Gr$  by
\begin{equation}\label{3.3}
h_t(g)={\tr}\;H_t(g),\quad g\in \Gr,\; t\ge0,
\end{equation}
where tr denotes the trace function on ${\End}(V_\sigma)$. Then it follows
from \ref{3.1} that $h_t$ belongs to ${\Co}^1(\Gr)$. In addition, 
$h_t$ is both left and right $K_\infty$-finite.

Let $\pi$ be a unitary irreducible representation of $\Gr$ on the Hilbert
space ${\mathcal H}(\pi)$. Then $\pi(h_t)$ is trace class and it follows as in 
\cite{BM} that
\begin{equation}\label{3.4}
\Tr\pi(h_t)=e^{t\pi(\Omega)}\dim 
\Hom_{K_\infty}({\cal H}(\pi),V_\sigma^*).
\end{equation}

\section[The main result]{The main result}
\setcounter{equation}{0}

First we note that by our assumptions on G,  strong approximation holds for 
$G(\A)$ and hence we have
$$L_{\cu}^2(G(\Q)\ba G(\A)/K_{\fin})\cong L^2_{\cu}(\Gamma\ba \Gr)$$
as $\Gr$-modules.
Let $h_t$, $t>0$, be the one-parameter
family of functions on $\Gr$ defined by (\ref{3.3}). By \cite[Corollary 0.2]
{Mu1} the restriction of $\rho_\infty(h_t)$ to $L^2_{\cu}(\Gamma\ba \Gr)$ 
is a trace class operator. However, since the support of $h_t$ is not compact
it is not known that we may insert $h_t$ in the trace formula for an 
arbitrary $G$ (although this is likely since
$h_t$ is a very rapidly decreasing function at infinity). 
This difficulty is bypassed by  modifying $h_t$ in the following way. Let 
$\varphi\in C^\infty_c(\R)$ be such that $\varphi(u)=1$, if $|u|\le1/2$, and
$\varphi(u)=0$, if $|u|\ge 1$. Given $t>0$, let 
$\varphi_t\in C_c^\infty(\Gr)$ be defined by
$$\varphi_t(g)=\varphi(r^2(g)/t^{1/2}),$$
where $r(g)$ is defined by (\ref{3.2}). Set
$$\tilde h_t(g)=\varphi_t(g)h_t(g),\quad g\in \Gr,\; t>0.$$
Then $\tilde h_t\in C^\infty_c(\Gr)$.  Recall that $\rho_\infty$ denotes
the right regular representation of $\Gr$ in $H_{\cu}^\Gamma$.

\begin{prop}\label{p4.1}
There exist $C,c>0$ such that
$$
\big|\Tr(\rho_\infty(h_t))-  
\Tr(\rho_\infty(\tilde h_t))\big|\le C e^{-c/\sqrt{t}}$$
for all $0< t\le 1$.
\end{prop}

\begin{proof}
For $\pi\in\Pi(\Gr)$ let $m_\Gamma(\pi)$ denote the multiplicity with 
which $\pi$ occurs in the discrete subspace $L^2_{\cu}(\Gamma\ba \Gr)$. 
For every $\pi\in\Pi(\Gr)$ with $m_\Gamma(\pi)\neq0$ we choose an 
orthonormal basis $e_{\pi,i}$, $i\in I(\pi)$, for the finite-dimensional
vector space $({\mathcal H}(\pi)\otimes V_\sigma)^{K_\infty}$. 
Let $\Delta_G$ be the Laplace operator on $\Gr$ with respect to the left
invariant metric. Note that 
$$\Delta_G=R(-\Omega+2\Omega_{K_\infty}),$$
where $\Omega_{K_\infty}$ is the Casimir element of $K_\infty$. Let $d=\dim X$
and let $\lambda_\sigma$ be the Casimir eigenvalue of $\sigma$. We note that
$$-\lambda_\pi+2\lambda_\sigma\ge0$$
for all $\pi\in\Pi(\Gr)$ with 
$(\H(\pi)\otimes V_\sigma)^{K_\infty}\neq\{0\}$ (cf. Lemma 2.6 of \cite{DH}). 
Let $\chi_\sigma$ be the character of $\sigma$. 
Then for every $f\in{\Co}^1(\Gr)$ with $f=f\ast{\overline \chi}_
\sigma$ we get

\begin{equation*}
\begin{split}
\Tr\rho_{\infty}(f) & =\sum_{\pi\in\Pi(\Gr)}
m_\Gamma(\pi)\sum_{i\in I(\pi)}(\rho_\infty(f)e_{\pi,i},e_{\pi,i})\\
& = \sum_{\pi\in\Pi(\Gr)}m_{\Gamma}(\pi)
(1-\lambda_\pi+2\lambda_\sigma)^{-4d}\\
& \quad\quad\times\sum_{i\in I(\pi)}\int_{\Gr}
\big (({\Id}+\Delta_G)^{4d}f\big )(g)(\rho_\infty(g)e_{\pi,i},e_{\pi,i})\;dg.\\
\end{split}
\end{equation*}

By \cite[Theorem 0.1]{Mu1} we have
$$\sum_{\pi\in\Pi(G(\R))}m_\Gamma(\pi)\dim(\H(\pi)\otimes V_\sigma)^{K_\infty}
(1-\lambda_\pi+2\lambda_\sigma)^{-4d}<\infty.$$
Using this result we get
\begin{equation}\label{4.1}
\big |\Tr\rho_{\infty}(f)\big |\le 
C\parallel({\Id}+\Delta_G)^{4d}f\parallel_{L^1(\Gr)}
\end{equation}

for some constant $C>0$, independent of $f$. Let 
$\psi_t=1-\varphi_t$, $t>0$. Then it follows from (\ref{4.1}) that
\begin{equation}\label{4.2}
\begin{split}
\big |\Tr(\rho_{\infty}(h_t))-\Tr(\rho_{\infty}(\tilde h_t))\big |
&=\big|\Tr(\rho_\infty(\psi_th_t)\big|\\
&\le C\parallel({\Id}+\Delta_G)^{4d}(\psi_t h_t)\parallel_{L^1(\Gr)}.
\end{split}
\end{equation}

It remains to estimate the right hand side. Let $X_1,...,X_a$ be an orthonormal
basis of $\gf(\R)$. Then $\Delta_G=-\sum_iX_i^2$. Denote by $\nabla$ the
canonical left invariant connection on $\Gr$. Then it follows that there
exists $C_1>0$ such that for all $f\in C^\infty(\Gr)$ we have
\begin{equation}\label{4.3}
\big |({\Id}+\Delta_G)^{4d}f(g)\big |\le C_1\sum_{j=0}^{8d}\parallel
\nabla^j f(g)\parallel,\quad g\in \Gr.
\end{equation}

By Proposition \ref{p3.1} there exist constants $C,c>0$ such that
$$\parallel\nabla^j h_t(g)\parallel\le C t^{-(a+j)/2}e^{-cr^2(g)/t},\quad
g\in \Gr,$$
for $j\le 8d$ and $0< t\le1$. Let $\chi_t$ be the characteristic function
of the set $\R-(-t^{1/4},t^{1/4})$. Recall that 
$\psi_t(g)=(1-\varphi)(r^2(g)/t^{1/2})$ and $(1-\varphi)(u)$ is constant for
$|u|\ge 1$. This implies that there exists a constant $C_2>0$ such that
\begin{equation}\label{4.4}
\parallel\nabla^j \psi_t(g)\parallel\le C_2 t^{-4d}\chi_t(r(g)),\quad 
g\in \Gr,
\end{equation}

for $j\le 8d$ and $0< t\le1$. Combining (\ref{4.3}) and (\ref{4.4}) we obtain

\begin{equation*}
\begin{split}
\sum_{j=0}^{8d}\parallel\nabla^j(\psi_th_t)(g)\parallel
& \le C_2 t^{-a/2-8d}\chi_t(r(g))e^{-c_2 r^2(g)/t}\\
& \le C_3e^{-c_3/\sqrt{t}} e^{-c_3 r^2(g)}\\
\end{split}
\end{equation*}

for all $g\in \Gr$ and $0<t\le1$. Finally note that for every $c>0$,
$e^{-c r^2(g)}$ is an integrable function on $\Gr$.
Together with (\ref{4.2}) and (\ref{4.3}) the proof follows.
\end{proof}

We are now ready to apply the trace formula. 
Given $\pi\in\Pi(\Gr)$, let $m_\Gamma(\pi,S)$ denote 
the multiplicity with which $\pi$ occurs in the regular representation of
$\Gr$ on $H^{\Gamma,S}_{\cu}$. Then it follows from (\ref{3.4}) that

\begin{equation}\label{4.5}
\Tr\left(\rho_\infty(h_t)|H^{\Gamma,S}_{\cu}\right)=\sum_{\pi\in\Pi(\Gr)}
e^{t\lambda_\pi}m_\Gamma(\pi,S)
\dim\Hom_{K_\infty}({\mathcal H}(\pi),V_\sigma^*).
\end{equation}

We note that the $\lambda_\pi$'s are the eigenvalues of the operator
 which is induced by $\rho_\infty(\Omega)$ in $H^{\Gamma,S}_{\cu}$ and 
$ m_\Gamma(\pi,S)\dim\Hom_{K_\infty}({\mathcal H}(\pi),V_\sigma^*)$ is the 
multiplicity of the corresponding eigenvalue.
Let 
$$N^\Gamma_{\cu}(T,\sigma,S)=\sum_{|\lambda_\pi|\le T}
 m_\Gamma(\pi,S)\dim\Hom_{K_\infty}({\mathcal H}(\pi),V_\sigma^*).$$
We may now state and prove our main result.
\begin{theo}\label{th4.2}  
Let $G$ be an almost simple, connected
and simply connected algebraic group defined over $\Q$ such that
$\Gr$ is non compact. Let $d_S$ denote the formal dimension
of the Steinberg representation of $G_S$ and let 
$$c_S(\Gamma)=c(K_S)d_S\vol(K_S).$$
Then, if $Z_\Gamma=\Gamma\cap Z_G(\R)\subset K_\infty$ is such that
$\sigma|_{Z_\Gamma}=
\Id$, we have
$$\lim_{T\to\infty}\frac{N^\Gamma_{cusp}(T,\sigma,S)}{T^{d/2}}= 
\;c_S(\Gamma)\frac{\vol(\Gamma\ba X)}{(4\pi)^{d/2}
{\bf \Gamma}(d/2+1)}.$$
\end{theo}

\begin{proof}
We apply the trace formula to
$$f_t=\tilde h_t\otimes e_{\fin,S}\otimes f_S.$$
The assumptions of Corollary \ref{c2.3} are satisfied.
We first observe that  $b(\tilde h_t,S)$ remains bounded. Since
$\Gr$ is non compact we have 
$$\dim (\Gr/K_\infty)=d>0$$ 
and hence
$t^{d/2}b(\tilde h_t,S)\to0$ as $t\to0$. Then Corollary \ref{c2.3} yields
$$t^{d/2}\tr(\rho_\infty(\tilde h_t)|H^{\Gamma,S}_{\cu})=
c(K_S)t^{d/2}\sum_{\gamma\in\mG_e} a(\gamma){\cO}_\gamma(f_t)+o(1)$$
as $t\to0$. Since the support of $f_\Sigma$ remains in a fixed compact set
we already observed in Proposition \ref{p2.1} that in
the sum over $\mG_e$ only a finite fixed number of terms may be
non zero. Hence it suffices to estimate each term. We first remark that
$${\cO}_\gamma(f_t)={\cO}_\gamma(\tilde h_t){\cO}_\gamma(e_{\fin,S})
{\cO}_\gamma(f_S).$$
We shall now use results that are established for the case $G=\SL(n)$ in 
\cite{Mu2}, 
but notation make sense and proofs extend readily for the general case.
According to \cite[Proposition 7.3]{Mu2}, we have
$$t^{d/2}{J}_M(\gamma,\tilde h_t)\to0,\quad t\to0,$$ 
for $M=G$ and $\gamma\notin Z_G(\R)$, 
while we know by \cite[Lemma 2.3]{Mu2} that $t^{d/2}h_t(z)$ has a non zero 
positive limit for $z\in Z_G(\R)$~:
\begin{equation}\label{4.6}
t^{d/2}h_t(z)\to\frac{\tr\sigma(z)}{(4\pi)^{d/2}\vol(K_\infty)}.
\end{equation}
Since 
$$J_G(\gamma,h)={\cal O}_\gamma(h)$$
 it follows that non central elements give a negligible contribution
 and, since $d_S=f_S(z)$~, we get 
$$ \lim_{t\to0}t^{d/2}\Tr(\rho_\infty(\tilde h_t)|H^{\Gamma,S}_{\cu})=
\frac{\vol(G(\Q)\ba G(\A))\,
c_S(\Gamma)}{\vol(K_{\fin})}
\lim_{t\to0} t^{d/2}\sum_{z\in Z_\Gamma}h_t(z)$$
 with 
$$ c_S(\Gamma)=c(K_S)\, d_S\vol(K_S)\,\,.$$
We observe that
$$\frac{\vol(G(\Q)\ba G(\A))}{\vol(K_\infty)\vol(K_{\fin})}
=\frac{\vol(\Gamma\ba X)}{\text{Card}(Z_\Gamma)}.$$
Using (\ref{4.6}) and the assumption that $\sigma|_{Z_\Gamma}=\Id$, it follows
that the above limit can be rewritten as
$$ \lim_{t\to0}t^{d/2}\tr(\rho_\infty(\tilde h_t)|H^{\Gamma,S}_{\cu})
=\,c_S(\Gamma)\frac{\dim\sigma\;\vol(\Gamma\ba X)}{(4\pi)^{d/2}}.$$
Now, it follows from Proposition \ref{p4.1} that
$$ \lim_{t\to0}t^{d/2}\Tr(\rho_\infty(\tilde h_t)|H^{\Gamma,S}_{\cu})=
\lim_{t\to0} t^{d/2}\Tr(\rho_\infty(h_t)|H^{\Gamma,S}_{\cu}),$$ 
and hence we get
$$\lim_{t\to0} t^{d/2}\Tr(\rho_\infty(h_t)|H^{\Gamma,S}_{\cu})
=\,c_S(\Gamma)\frac{\dim\sigma \vol(\Gamma\ba X)}{(4\pi)^{d/2}}.$$
We conclude using (\ref{4.5}) and the Tauberian theorem.
\end{proof}

As explained in the introduction, Theorem \ref{th1.1} follows immediately 
from Theorem \ref{th4.2}.

{\bf Remarks.} 1) The constant $c_S(\Gamma)$ is non zero if and only if
$c(K_S)\ne0$ and this is true whenever $K_S$ is small enough so that
the space of the Steinberg representation of $G_S$
contains non zero $K_S$-invariant vectors. This is the case
if $K_p\subset I_p$ a minimal parahoric subgroup for $p\in S$.

2) If the assumption $\sigma|_{Z_\Gamma}=\Id$ in Theorem \ref{4.2} is not
satisfied, i.e, if there exists $z\in Z_\Gamma$ such that $\sigma(z)\not=\Id$,
then $H^\Gamma_{\cu}(\sigma)=\{0\}$.

\end{document}